\documentclass[12pt]{amsart}
\usepackage{amsmath,amscd,amssymb,amsfonts}
\newcommand{\h}{\hbox}
\newcommand{\q}{\quad}

\newcommand{\nin}{\par\noindent}
\newcommand{\bs}{\par\bigskip}
\newcommand{\ms}{\par\medskip}
\newcommand{\sk}{\par\smallskip}
\newcommand{\mopl}{\h{$\bigoplus$}}
\newcommand{\msum}{\h{$\sum$}}
\newcommand{\mprod}{\h{$\prod$}}
\newcommand{\mcap}{\h{$\bigcap$}}
\newcommand{\mcup}{\h{$\bigcup$}}
\newcommand{\al}{\alpha}
\newcommand{\C}{{\mathbf C}}
\newcommand{\Cc}{{\mathcal C}}

\newcommand{\De}{\Delta}
\newcommand{\Dt}{\widetilde{D}}
\newcommand{\Det}{\widetilde{\Delta}}
\newcommand{\E}{\widetilde{E}}

\newcommand{\ft}{\tilde{f}}
\newcommand{\Ga}{\Gamma}
\newcommand{\ga}{\gamma}
\newcommand{\HH}{{\mathcal H}}
\newcommand{\iti}{\tilde{i}}
\newcommand{\jti}{\tilde{j}}

\newcommand{\la}{\lambda}
\newcommand{\lab}{\overline{\lambda}}
\newcommand{\M}{{\mathcal M}}

\newcommand{\OO}{{\mathcal O}}
\newcommand{\om}{\omega}

\newcommand{\PP}{{\mathbf P}}
\newcommand{\Q}{{\mathbf Q}}
\newcommand{\R}{{\mathbf R}}
\newcommand{\Tt}{\widetilde{T}}
\newcommand{\Xt}{\widetilde{X}}
\newcommand{\Yt}{\widetilde{Y}}
\newcommand{\Zt}{\widetilde{Z}}
\newcommand{\Z}{{\mathbf Z}}
\newcommand{\Gr}{{\rm Gr}}
\newcommand{\IH}{{\rm IH}}
\newcommand{\Ker}{{\rm Ker}}
\newcommand{\Sp}{{\rm Sp}}
\newcommand{\into}{\hookrightarrow}
\newcommand{\bl}{\bigl}
\newcommand{\br}{\bigr}
\newcommand{\da}{\raise-2pt\h{$\downarrow$}}
\newcommand{\ssb}{\raise.15ex\h{${\scriptscriptstyle\bullet}$}}
\newcommand{\ssc}{\,\raise.15ex\h{${\scriptstyle\circ}$}\,}
\newcommand{\simto}{\buildrel{\sim}\over\longrightarrow}

\begin{document}
\title[Number of Jordan blocks of the maximal size]
{Number of Jordan blocks of the maximal size\\
for local monodromies}
\author{Alexandru Dimca}
\address{Institut Universitaire de France et
Laboratoire J.A.\ Dieudonn\'e, UMR du CNRS 7351,
Universit\'e de Nice-Sophia Antipolis, Parc Valrose,
06108 Nice Cedex 02, France}
\email{Alexandru.DIMCA@unice.fr}
\author{Morihiko Saito}
\address{RIMS Kyoto University, Kyoto 606-8502 Japan}
\email{msaito@kurims.kyoto-u.ac.jp}
\dedicatory{Dedicated to Professor Joseph Steenbrink}
\begin{abstract}
We prove formulas for the number of Jordan blocks of the maximal size
for local monodromies of one-parameter degenerations of complex
algebraic varieties where the bound of the size comes from the
monodromy theorem.
In case the general fibers are smooth and compact, the proof calculates
some part of the weight spectral sequence of the limit mixed Hodge
structure of Steenbrink.
In the singular case, we can prove a similar formula for the monodromy
on the cohomology with compact supports, but not on the usual cohomology.
We also show that the number can really depend on the position of
singular points in the embedded resolution even in the isolated
singularity case, and hence there are no simple combinatorial
formulas using the embedded resolution in general.
\end{abstract}
\maketitle

\centerline{\bf Introduction}

\bs\nin
Let $f:X\to\De$ be a proper surjective morphism of a connected complex
manifold $X$ to an open disk $\De$, which is smooth over $\De^*$.
Assume there is a proper surjective morphism from a K\"ahler manifold
to $X$. One may assume for simplicity that $f$ is a projective
morphism. For a divisor $D$ on $X$, set
$$U:=X\setminus D,\q f_U:=f|_U:U\to\De,\q U_t:=f_U^{-1}(t),\q
X_t:=f^{-1}(t).$$
Shrinking $\De$ if necessary, we may assume that the
$H^j(U_t,\Q)\,\,(t\in\De^*)$ form local systems and moreover
$H^j(U_t,\Q)=(R^j(f_U)_*\Q_U)_t\,\,(t\in\De^*)$ for any $j$.
We are interested in their monodromy around the origin.
So we may assume that $X_0\cup D$ is a divisor with simple normal
crossings, and every irreducible components $D_k$ of $D$ is dominant
over $\De$. Then, shrinking $\De$ if necessary, we may assume moreover
that any $D_k$ and any intersections of $D_k$ are smooth over $\De^*$.
Let $Y_i$ be the irreducible components of $Y:=X_0\subset X$ with $m_i$
the multiplicity of $Y$ at the generic point of $Y_i$.
Set $Y_I:=\mcap_{i\in I}Y_i$.
\sk
Set $J(\la):=\{i\mid\la^{m_i}=1\}$ for a root of unity $\la$ in $\C^*$.
For $I\subset J(\la)$, let $Y^{(\la)}_I\subset Y_I$ be the union of
the connected components of $Y_I$ which do not intersect
$Y_{i'}$ for any $i'\notin J(\la)$.
Note that $Y^{(1)}_I=Y_I$ for $\la=1$.
We have the complex $C_{f,\la}^{\ssb}$ defined by
$$C_{f,\la}^j:=\mopl_{I\subset J(\la),|I|=j+1}\,H^0(Y^{(\la)}_I,\C),$$
where the differential is induced by the Cech restriction
morphism as in [St1].
Similarly we have $Y_{k,I}$, $Y^{(\la)}_{k,I}$, $C_{f_k,\la}^{\ssb}$
for each $k$ by replacing $f:X\to\De$ with $f_k:=f|_{D_k}:D_k\to\De$
and $Y_i$ with $Y_{k,i}:=D_k\cap Y_i$.
There are canonical restriction morphisms
$$r_k:C_{f,\la}^{\ssb}\to C_{f_k,\la}^{\ssb}.$$
Set $n:=\dim X-1$.
Let $\nu_{f_U,\la}^j$ (resp. $\nu_{c,f_U,\la}^j$) denote the number of
Jordan blocks of the theoretically maximal size $j+1$ for the
eigenvalue $\la$ of the monodromy on $H^j(U_t,\C)$
(resp. $H_c^j(U_t,\C)$) for $t\in\De^*$ and $j\in[0,n]$.
Here the upper bounds come from the monodromy theorem.
This upper bound is $2n-j$ for $j>n$, and the number of Jordan blocks of
the maximal size for the eigenvalue $\la$ of the monodromy on
$H^j(U_t,\C)$ and $H_c^j(U_t,\C)$ are respectively given by
$\nu_{c,f_U,\lab}^{2n-j}$ and $\nu_{f_U,\lab}^{2n-j}$ for $j\in[n,2n]$
by duality. Thus it is enough to consider $\nu_{f_U,\la}^j$,
$\nu_{c,f_U,\la}^j$ for $j\in[0,n]$ in the smooth case.
We have the following:
\ms\nin
{\bf Question~1.} Do the following equalities hold for $j\in[0,n]$?
$$\nu_{f_U,\la}^j=\dim H^jC_{f,\la}^{\ssb},\q
\nu_{c,f_U,\la}^j=\dim\Ker\bl(H^jC_{f,\la}^{\ssb}\to
\mopl_k\,H^jC_{f_k,\la}^{\ssb}\br).$$

These equalities follow from the theory of limit mixed Hodge structures
([St1],[StZ]) if $\la=1$, and they were expected to hold also for
$\la\ne 1$ if $f$ is obtained by a desingularization of a good
compactification of a germ of a holomorphic function $g_0$ with an
isolated singularity as in Theorem~3 below where $D=\emptyset$, i.e.
$f_U=f$. In fact, we can prove the equality
$\nu_{f,\la}^n=\nu_{g_0,\la}^n=(-1)^n\chi(C_{f,\la}^{\ssb})$ for
$\la\ne 1$ as in Theorem~3 below, for instance, in case of super-isolated
singularities [Lu], or more generally, Yomdin singularities [Yo] with
$n=2$, see Proposition~(3.8) below (and also [Ar1], [MM]).
However, it turns out that the answer to Question~$1$ is negative,
and there is no simple combinatorial formula in general,
since we have quite recently found the following:
\ms\nin
{\bf Theorem~1.} {\it 
The $\nu_{f,\la}^j$ cannot be determined only by the combinatorial data of
the embedded resolution, and may really depend on the position of the
singular points in the embedded resolution even in case $f$ is obtained by
a desingularization of a good compactification of a germ of a
holomorphic function with an isolated singularity.}
\ms
Here a good compactification means a compactification having only one
singular point as constructed in [Br], and the combinatorial data of
an embedded resolution in this paper mean the intersection lattice
consisting of the connected components of the $Y^{(\la)}_I$ with
$\la$ fixed (see also [Ar1]). Theorem~1 will be shown in (4.3) below.
In Theorem~3, it will be shown that the $\nu_{f,\la}^j$ are
determined by the dimensions of the $C_{f,\la}^j$ (i.e. by the numbers
of the connected components of the $Y^{(\la)}_I$ with $|I|=j+1$)
in case of a desingularization of a good compactification of a germ
of a holomorphic function with an isolated singularity, provided that
$B_{f,\la}^j=C_{f,\la}^j$ for any $j$ in the notation of
Theorem~2 below.
A sufficient condition for the last equality is given in Theorem~4(i).
Note that Theorem~1 is related with certain earlier work in the literature
like [Ar1], [Ar2], [AC], [GaNe], [GLM], [Li], [MH], [Za], etc.
\sk
We now explain an improvement of the above formula in Question~1.
Let $H^j(U_{\infty})_{\la}$ (resp. $H_c^j(U_{\infty})_{\la}$)
denote the $\la$-eigenspace of the limit mixed Hodge structure with
$\C$-coefficient. Let $W$ be the weight filtration of the limit mixed
Hodge structure. We have the following.
\ms\nin
{\bf Theorem~2.} {\it There are complexes $B_{f,\la}^{\ssb}$,
$B_{f_k,\la}^{\ssb}$ and morphisms
$r'_k:B_{f,\la}^{\ssb}\to B_{f_k,\la}^{\ssb}$ with
$B_{f,\la}^j$, $B_{f_k,\la}^j$ direct factors of $C_{f,\la}^j$,
$C_{f_k,\la}^j$ respectively, and such that we have for $j\in[0,n]$
$$\Gr^W_0H^j(U_{\infty})_{\la}=H^jB_{f,\la}^{\ssb},\q
\Gr^W_0H_c^j(U_{\infty})_{\la}=\Ker\bl(H^jB_{f,\la}^{\ssb}\to
\mopl_k\,H^jB_{f_k,\la}^{\ssb}\br),
\leqno(0.1)$$
\vskip-20pt
$$\nu_{f_U,\la}^j=\dim H^jB_{f,\la}^{\ssb},\q
\nu_{c,f_U,\la}^j=\dim\Ker\bl(H^jB_{f,\la}^{\ssb}\to
\mopl_k\,H^jB_{f_k,\la}^{\ssb}\br).
\leqno(0.2)$$
The differentials of $B_{f,\la}^{\ssb}$, $B_{f_k,\la}^{\ssb}$ are
induced by the Cech restriction morphisms up to some nonzero constant
multiples which may depend on each inclusion of connected components
with codimension $1$.
We have $B_{f,1}^{\ssb}=C_{f,1}^{\ssb}$ and
$B_{f_k,1}^{\ssb}=C_{f_k,1}^{\ssb}$ if $\la=1$.}
\ms
Here the problem is the global triviality of the local systems
$L_{\la,I}$ of rank 1 in (1.1.4) below which are associated with the
nearby cycles, and we get $B_{f,\la}^j$ by replacing $Y^{(\la)}_I$
in the definition of $C_{f,\la}^j$ with a union of the connected
components of $Y^{(\la)}_I$ on which $L_{\la,I}$ is trivial,
and choosing a trivialization of $L_{\la,I}$ (and similar for
$B_{f_k,\la}^j$).
\sk
In the proper case (i.e.\ $D=\emptyset$), Theorem~2 follows from
Steenbrink's construction of the limit mixed Hodge structures using
$V$-manifolds [St2] together with the theory of bi-graded modules of
Lefschetz-type [Sa1], Sect.~4 (see also [GuNa]).
Here we do not need [Sa1], 4.2.3.1 (i.e.\ [SaZ], 1.3.8),
since we use the lowest weight part of the $E_1$-complex
where only the Cech restriction morphisms appear.
The non-proper case then follows by using the limit of weight
spectral sequences in [StZ].
In Theorem~(2.2) below. Theorem~2 for $\nu_{c,f_U,\la}^j$ will be
generalized to the singular case although the assertion for
$\nu_{f_U,\la}^j$ cannot, see Example~(2.3) below.
\sk
It is not easy to determine the differential of the complex
$B_{f,\la}^{\ssb}$ for $\la\ne 1$ in general.
This problem can be avoided in the case of good compactifications
of isolated singularities as follows.

\ms\nin
{\bf Theorem~3.} {\it Assume $f$ is obtained by an embedded resolution
of a good compactification $g:X\to\De$ of a germ of a holomorphic
function $g_0:(\C^{n+1},0)\to(\C,0)$ with an isolated singularity.
Define $\nu_{g_0,\la}^n$ by using the Milnor monodromy where the
maximal size of Jordan blocks for $\la=1$ is $n$ instead of $n+1$.
Then we have for any $\la$
$$\aligned\nu_{f,\la}^n&=\nu_{g_0,\la}^n=
(-1)^n\bl(\chi(B_{f,\la}^{\ssb})-\delta_{\la,1}\br),\\
\nu_{f,\la}^j&=\delta_{\la,1}\,\delta_{j,0}\q(j\in[0,n-1]).
\endaligned$$
where $B_{f,1}^{\ssb}=C_{f,1}^{\ssb}$ in case $\la=1$.
Here $\chi(B_{f,\la}^{\ssb})$ is the Euler characteristic of the
complex $B_{f,\la}^{\ssb}$, and $\delta_{\al,\beta}$ is $1$ if
$\al=\beta$, and $0$ otherwise.}

\ms
It is quite difficult to determine $B_{f,\la}^j$, $B_{f_k,\la}^j$ for
$\la\ne 1$ in general. In fact, we may have
$\chi(B_{f,\la}^{\ssb})\ne\chi(C_{f,\la}^{\ssb})$ for $\la \ne 1$, and
moreover the inequality
$\dim H^jB_{f,\la}^{\ssb}\le\dim H^jC_{f,\la}^{\ssb}$ does not
necessarily hold, see Example~(4.1) and (4.3) below. The following
sufficient conditions for the coincidence are quite useful in certain
cases.
\ms\nin
{\bf Theorem~4.} {\it For $\la\ne 1$, set
$Y^{(\la)}=\mcup_{I\subset J(\la)}\,Y^{(\la)}_I$.
Let $m_{\la}$ be the order of $\la$.
\sk\nin
{\rm (i)} If $H^1(Y^{(\la)}_I,\Z/m_{\la}\Z)=0$ for any
$I\subset J(\la)$ with $|I|=j+1$, then $B_{f,\la}^j=C_{f,\la}^j$.
\sk\nin
{\rm (ii)} If $H^1(Y^{(\la)},\Z/m_{\la}\Z)=0$, then
$B_{f,\la}^{\ssb}=C_{f,\la}^{\ssb}$ as a complex.
\sk\nin
Similar assertions hold for $B_{f_k,\la}^{\ssb}$, $C_{f_k,\la}^{\ssb}$
by replacing respectively $Y^{(\la)}_I$, $Y^{(\la)}$ with
$Y^{(\la)}_{k,I}$,
$Y_k^{(\la)}:=\mcup_{I\subset J(\la)}\,Y^{(\la)}_{k,I}$.}
\ms
In (i) the problem is the triviality of certain unramified cyclic
covering of $Y^{(\la)}_I$ with degree $m_{\la}$, and hence the
cohomology $H^1(Y^{(\la)}_I,\Z/m_{\la}\Z)$ appears.
This may be replaced with $H_1(Y^{(\la)}_I,\Z)$ (since the
monodromy group of a local system of rank 1 is abelian), but not with
$H^1(Y^{(\la)}_I,\Z)$.
Similar assertions hold for (ii) with $Y^{(\la)}_I$ replaced by
$Y^{(\la)}$.

In the case of Theorem~3, Theorem~4(i) is enough since we do not
have to consider the differential in this case.
The hypothesis of Theorem~4(ii) is rather strong, and is not
often satisfied except for certain special cases, e.g. if $f$ is as in
Theorem~3 with $n=1$ (i.e. $\dim X=2$), and the embedded resolution
is obtained by repeating point-center blow-ups.
In this case, Theorem~3 for $\la\ne 1$ means
$$\aligned\nu_{g_0,\la}^1&=\#\bl\{I\subset J(\la)\,\big|\,|I|=2,\,
Y_I\ne\emptyset\br\}\\ &\q-\#\bl\{j\in J(\la)\,\big|\,Y_j\cap Y_i
=\emptyset\,\,\,\h{for any}\,\,i\notin J(\la)\br\}.\endaligned
\leqno(0.3)$$
For $\la=1$ and $n=1$, Theorem~3 simply gives a well-known
formula $\nu_{g_0,1}^1=r_{g_0}-1$ where $r_{g_0}$ is the number
of analytic local irreducible components of $g_0^{-1}(0)$.
\sk
Theorem~2 improves a result of Y.~Matui and K.~Takeuchi [MT]
where the number is bounded by $\dim C_{f,\la}^j$ in the case of
monodromies at infinity of polynomial maps with $\la\ne 1$
(since $\dim C_{f,\la}^j\ge\dim H^jB_{f,\la}^{\ssb}$).
In case $j=n$, the latter assertion easily follows from a
{\it local} assertion at the level of perverse sheaves in [Sa5], 3.2.2:
$$\min\bl\{k>0\,\big|\,N^k(\psi_{f,\la}\C_X)=0\,\,\,
\h{around}\,\,\,x\br\}=\#\bl\{i\in J(\la)\mid x\in D_i\br\},
\leqno(0.4)$$
where $N$ is the nilpotent part of the monodromy $T$, and
$\psi_{f,\la}\C_X$ is the $\la$-eigenspace of the nearby cycles
$\psi_f\C_X$ which is a shifted perverse sheaf,
see also [DS], 1.4 for a more precise local structure. This is more
or less well-known to the specialists of limit mixed Hodge structures
who are familiar with the theory of Steenbrink in [St2].
A more precise local structure as in [DS], 1.4 is implicit in the
definition of motivic Milnor fibers, and was used in the proof of the
compatibility with the Hodge realization by Denef and Loeser [DL].
\sk
The rank of the differential of $B_{f,\la}^{\ssb}$
as well as the difference between $\dim H^jB_{f,\la}^{\ssb}$ and
$\dim C_{f,\la}^j$ can be quite large as is seen in the case of
Example~(4.4) below.
This example shows that, even in the non-degenerate Newton
boundary case, we have to apply many blow-ups in order to get a
divisor with normal crossings (in the usual sense) by taking a
suitable subdivision of the dual fan, and the estimate in [MT] may
become rather bad (unless the dual fan is already smooth,
i.e.\ consisting of simplicials generated by integral vectors with
determinant 1).
The situation seems to be similar in the case of monodromies at
infinity.

\sk
In Section~1 we recall some basics of nearby cycles and limit mixed
Hodge structures in the non-reduced case, and then prove Theorems~2
and 3.
In Section~2 we partially generalize Theorem~2 to the singular case
in Theorem~(2.2).
In Section~3 we provide a method to show Theorem~1 in Section~4,
and prove Theorem~4.
In Section~4 we give some interesting examples, and prove Theorem~1
in (4.3).

\sk
The first named author was partially supported by the grant
ANR-08-BLAN-0317-02 (SEDIGA).
The second named author was partially supported by Kakenhi 21540037.

We thank the referee for useful comments especially about the
references.
\bs\bs
\centerline{\bf 1. Nearby cycles and limit mixed Hodge structures}

\bs\nin
In this section we recall some basics of nearby cycles and limit mixed
Hodge structures in the non-reduced case, and then prove Theorems~2
and 3.
\ms\nin
{\bf 1.1.~Local structure of nearby cycle sheaves.}
Let $f$ be a nonconstant holomorphic function on a complex manifold
$X$ of dimension $n+1$. Let $\psi_f\C_X$ denote the nearby cycle
sheaf with monodromy $T$ in [De2].
It is well known that this is a shifted perverse sheaf [BBD] (i.e.
$\psi_f\C_X[n]$ is a perverse sheaf).
Using the minimal polynomial of $T$, we have the Jordan decomposition
$T=T_sT_u$, where $T_s$ and $T_u$ respectively denote the semisimple
and unipotent part. For $\la\in\C^*$, set
$$\psi_{f,\la}\C_X:=\Ker(T_s-\la)\subset\psi_f\C_X,$$
in the abelian category of {\it shifted} perverse sheaves [BBD].
Then $\psi_{f,\la}\C_X=0$ except for a finite number of $\la$ which
are roots of unity, and
$$\psi_f\C_X=\mopl_{\la}\,\psi_{f,\la}\C_X.$$
Set $N=(2\pi i)^{-1}\log T_u$.
The weight filtration $W$ on $\psi_f\C_X$ is given by the monodromy
filtration with center $0$, i.e.
$$N^k:\Gr^W_k\psi_f\C_X\simto\bl(\Gr^W_{-k}\psi_f\C_X\br)(-k)\q(k>0),
\leqno(1.1.1)$$
where $(-k)$ is the Tate twist which shifts the weights by $2k$.
Define the $N$-primitive part by
$$P\Gr_k^{W}\psi_{f,\la}\C_X:=\Ker\,N^{k+1}\subset
\Gr_k^W\psi_{f,\la}\C_X\q(k\ge 0),$$
where it is zero for $k<0$.
By (1.1.1) we have the primitive decomposition
$$\Gr_j^W\psi_{f,\la}\C_X =
\mopl_{k\ge 0}\,N^k\bl(P\Gr_{j+2k}^W\psi_{f,\la}\C_X\br)(k).
\leqno(1.1.2)$$

Assume that $Y:=f^{-1}(0)$ is a divisor with simple normal crossings.
Let $Y_i$ be the irreducible components of $Y$ with multiplicities
$m_i$. Set $Y_I:=\mcap_{i\in I}Y_i$. For a root of unity $\la$ in
$\C^*$, set
$$J(\la):=\{i\mid\la^{m_i}=1\}.
\leqno(1.1.3)$$
By [Sa2], 3.3 (see also [DS], 1.4) we have the decomposition
$$P\Gr_k^{W}\psi_{f,\la}\C_X=\mopl_{I\subset J(\la),|I|=k+1}\,
(j_{\la,I})_{!*}L_{\la,I}(-k)[n-k].
\leqno(1.1.4)$$
Here $L_{\la,I}$ is a local system of rank $1$ underlying a locally
constant variation of complex Hodge structure of weight $0$ on
$U_{\la,I}:=Y_I\setminus\mcup_{i\notin J(\la)}\,Y_i$, and
$(j_{\la,I})_{!*}$ is the intermediate direct image [BBD] by the
natural inclusion $j_{\la,I}:U_{\la,I}\into Y_I$.
Furthermore, the monodromy of $L_{\la,I}$ around
$Y_{j}\,(j \notin J(\la))$ is given by the multiplication by
$\la^{-m_j}$ so that
$$(j_{\la,I})_{!*}L_{\la,I}[n-k]=
(j_{\la,I})_!L_{\la,I}[n-k]=
\R(j_{\la,I})_*L_{\la,I}[n-k].
\leqno(1.1.5)$$
Indeed, the last isomorphism follows from the above
information of the local monodromies, and the first
isomorphism follows from this by the definition of
the intermediate direct image $(j_{\la,I})_{!*}$
(see [BBD]).

\ms\nin
{\bf 1.2.~Relation with Steenbrink's construction.}
In the above notation and assumption, let $\Xt$ be the normalization
of the base change of $f:X\to\De$ by the totally ramified $m$-fold
covering $\Det\to\De$ with $m:={\rm LCM}(m_i)$, see [St2].
Let $\pi:\Xt\to X$, $\ft:\Xt\to\Det$ be the canonical morphisms.
Set $\Yt:=\pi^{-1}(Y)$, and let $\pi_0:\Yt\to Y$ be the restriction
of $\pi$ over $Y$. Then we have a canonical isomorphism
$$\psi_f\C_X=(\pi_0)_*\psi_{\ft}\C_{\Xt},
\leqno(1.2.1)$$
where the monodromy $\Tt$ on the right-hand side is identified with
the $m$-th power of the monodromy $T$ on the left-hand side
which is unipotent.
This follows from the commutative diagram
$$\begin{matrix}\Yt&\buildrel{\iti}\over\into&\Xt&
\buildrel{\,\,\jti_{\infty}}\over\longleftarrow&\Xt_{\infty}\\
\,\,\,\,\da{\scriptstyle\pi_0}&&
\,\,\,\da{\scriptstyle\pi}&&
\,\,\,\da{\scriptstyle\pi_{\infty}}\\
Y&\buildrel{i}\over\into&X&
\buildrel{\,\,j_{\infty}}\over\longleftarrow&X_{\infty}\end{matrix}
\leqno(1.2.2)$$
where $X_{\infty}$ is the base change of $X$ by the universal covering
of $\De^*$ over $\De$, and similarly for $\Xt_{\infty}$ with $X$,
$\De^*$, $\De$ replaced by $\Xt$, $\Det^*$, $\Det$.
Here $\pi_{\infty}$ is an isomorphism.
Then (1.2.1) follows from the definition of the nearby cycles
$\psi_f\C_X:=i^*\R(j_{\infty})_*\C_{X_{\infty}}$ (and similarly for
$\psi_{\ft}\C_{\Xt}$) by using the diagram (1.2.2) together with
the commutativity
$$(\pi_0)_*\ssc\iti^*=i^*\ssc\pi_*,$$
where $\pi$ is finite and hence proper.
The relation between $\Tt$ and $T^m$ is clear by the construction
of the isomorphism.
\sk
By the above construction, the Milnor fiber of $\ft$ at any point
of $\Yt$ is connected, and we have
$$\HH^0\psi_{\ft}\C_{\Xt}=\C_{\Yt}.$$
Combining this with (1.2.1), we get
$$\HH^0\psi_f\C_X=(\pi_0)_*\C_{\Yt},$$
since $\pi_0$ is finite. The action of $T$ on the left-hand side
is semisimple and corresponds to the action of an appropriate
generator of the covering transformation group $\Z/m\Z$ of $\pi$.
So we get
$$\HH^0\psi_{f,\la}\C_X=\bl((\pi_0)_*\C_{\Yt}\br)_{\la},
\leqno(1.2.3)$$
where the right-hand side denotes the $\la$-eigenspace.
\sk
Set $\Yt_I:=\pi^{-1}(Y_I)$. We have the Cech resolution
$$\C_{\Yt}\simto\Cc_{\Yt}^{\ssb}\q\h{with}\q
\Cc_{\Yt}^j:=\mopl_{|I|=j+1}\,\C_{\Yt_I}.
\leqno(1.2.4)$$
Taking the direct image by $(\pi_0)_*$ and the $\la$-eigenspace,
we get the quasi-isomorphism
$$\bl((\pi_0)_*\C_{\Yt}\br)_{\la}\simto\Cc_{Y,\la}^{\ssb}:=
\bl((\pi_0)_*\Cc_{\Yt}^{\ssb}\br)_{\la}.
\leqno(1.2.5)$$
In the notation (1.1.4), we have moreover
$$\Cc_{Y,\la}^j=\mopl_{I\in J(\la),|I|=j+1}\,(j_{\la,I})_{!*}L_{\la,I},
\leqno(1.2.6)$$
since
$$\bl((\pi_0)_*\C_{\Yt_I}\br)_{\la}=\begin{cases}(j_{\la,I})_{!*}L_{\la,I}&
\h{if}\,\,\,I\in J(\la),\\ \,0&\h{if}\,\,\,I\notin J(\la).\end{cases}
\leqno(1.2.7)$$
This can be reduced to the case $|I|=1$ by choosing any $i\in I$
and using the restriction morphisms
$$(j_{\la,I'})_{!*}L_{\la,I'}\to(j_{\la,I})_{!*}L_{\la,I}\q\h{for}\,\,\,
I'\subset I\in J(\la).
\leqno(1.2.8)$$
Moreover, we may restrict to any dense Zariski-open subset of
$Y_I$ (e.g.\ to the smooth points of $Y_{\rm red}$ if $|I|=1$) by
using the intermediate direct image by the open inclusion of the
Zariski-open subset, since the intermediate direct images commute
with the direct image by any finite morphisms (e.g. $(\pi_0)_*$),
see [BBD]. Then the assertion follows from (1.1.4) and (1.2.3).

\ms\nin
{\bf 1.3.~Weight spectral sequences.}
With the notation of (1.1) assume $f:X\to\De$ is a projective morphism
to an open disk $\De$, and moreover $Y$ is a divisor with simple
normal crossings.
Then we have the weight spectral sequence
$$E_1^{-k,j+k}=H^j\bl(Y,\Gr^W_k\psi_{f,\la}\C_X\br)
\Longrightarrow H^j(X_{\infty})_{\la},
\leqno(1.3.1)$$
where $W$ on $\psi_{f,\la}\C_X$ is the monodromy filtration as in
(1.1), and $H^j(X_{\infty})_{\la}$ denotes the $\la$-eigenspace of
the limit mixed Hodge structure with complex coefficients
as in [St1], [St2].
By [Sa1], Section~4 (or [GuNa]) the filtration on $H^j(X_{\infty})$
induced by $W$ on $\psi_{f,\la}\C_X$ is also the monodromy filtration
with center $0$, and we need the {\it shift} by $j$ to get the
weight filtration $W$ on $H^j(X_{\infty})$.

From the primitive decomposition (1.1.2) together with (1.1.4),
we can deduce the double complex structure of the $E_1$-complex in
[Sa1], Sect.~4 (see also [SaZ], 1.1) as follows:
$$\aligned&\,E_1^{-k,j+k}=\mopl_{a-b=k}\,C_{\la,a,b}^j\q\q\h{with}\\
&C_{\la,a,b}^j:=\begin{cases}\mopl_{|I|=a+b+1}\,\IH^{j-a-b}(Y_I,L_{\la,I})
(-a)&\h{if}\,\,\,a,b\ge 0,\\\,\,0&\h{otherwise},\end{cases}\endaligned
\leqno(1.3.2)$$
where $\IH^{j-a-b}(Y_I,L_{\la,I})$ is the intersection cohomology
[BBD], the action of $N$ is induced by
$$id:C_{\la,a,b}^j\to C_{\la,a-1,b+1}^j(-1)\q\h{if}\,\,\,
a-1,b\in\Z_{\ge 0},
\leqno(1.3.3)$$
and the $E_1$-differential is the sum of
$$d':C_{\la,a,b}^j\to C_{\la,a-1,b}^{j+1}\q\h{and}\q
d'':C_{\la,a,b}^j\to C_{\la,a,b+1}^{j+1},
\leqno(1.3.4)$$
which are identified up to a certain sign with the morphisms induced
respectively by the Cech-Gysin morphisms $\widetilde{\ga}$ and the
Cech restriction morphisms $\widetilde{\rho}$ between the $\Yt_I$
using the isomorphism (1.2.7).
In other words, (1.3.2) is obtained by taking the $\la$-eigenspace
of the direct image by $\pi_0$ of the double complex structure of
the $E_1$-complex for $\ft$ in [St2].
(Note that the kernel and image filtrations $K_i$ and $I^k$ in [SaZ]
are defined respectively by the conditions $a\le i$ and $b\ge k$.)
\sk
Consider now the {\it lowest} weight part of the $E_1$-complex.
Its weight $j+k$ is zero with
$$a=j-b=0\q\h{in}\,\,\,(1.3.2),$$
since the weight of $\IH^{j-a-b}(Y_I,L_{\la,I})(-a)$ is $(j-a-b)+2a$
with
$$j-a-b\ge 0,\q a\ge 0.$$
So the lowest weight part is the complex with $j$-th component given by
$$C_{\la,0,j}^j=\HH^0(Y,\Cc_{Y,\la}^j),
\leqno(1.3.5)$$
where the last isomorphism comes from (1.2.6).

\ms\nin
{\bf 1.4.~Limits of weight spectral sequences.}
With the notation of (1.1), let $D$ be a divisor with simple normal
crossings on $X$ such that all the irreducible components $D_j$ of
$D$ are dominant over $\De$. Set $U:=X\setminus D$, and
$D_J:=\mcap_{j\in J}\,D_j$ (where $D_{\emptyset}=X$).
We have the spectral sequences of mixed $\Q$-Hodge structures
compatible with the action of the semisimple part $T_s$ of the
monodromy:
$$\aligned{}_{\infty}E_1^{-i,j+i}&=\mopl_{|J|=i}\,H^{j-i}
(D_{J,\infty})(-i)\Longrightarrow H^j(U_{\infty}),\\
{}_{\infty}\!{}^cE_1^{i,j-i}&=\mopl_{|J|=i}\,H^{j-i}(D_{J,\infty})
\Longrightarrow H^j_c(U_{\infty}),\endaligned
\leqno(1.4.1)$$
which are dual of each other. They degenerate at $E_2$ since they are
the `limit' by $t\to 0$ of the weight spectral sequences
$$\aligned{}_tE_1^{-i,j+i}&=\mopl_{|J|=i}\,H^{j-i}(D_{J,t})(-i)
\Longrightarrow H^j(U_{t}),\\
{}_t\!{}^cE_1^{i,j-i}&=\mopl_{|J|=i}\,H^{j-i}(D_{J,t})
\Longrightarrow H^j_c(U_{t}),\endaligned
\leqno(1.4.2)$$
where the nearby cycle functor $\psi$ of mixed Hodge modules can be used
to define the `limit'. Here $D_{J,t}:=D_J\cap X_{t}$ for $t\in\De^*$.
The first spectral sequence in (1.4.1) was obtained in [StZ] in the
unipotent monodromy case, and it can be generalized to the
non-unipotent case by [St2].
Here the `limit' can be defined also by using the nearby cycle functor
$\psi$ of mixed Hodge modules and the spectral sequences are defined
by the weight filtration on the shifted perverse sheaves $(j_U)_*\Q_U$
or $(j_U)_!\Q_U$ with $j_U:U\into X$ the natural inclusion.
This implies for instance
$$\Gr^W_0(j_U)_*\Q_U=\Gr^W_0(j_U)_!\Q_U=\Q_X.$$
The $E_1$-differential of the spectral sequences are induced
by the Cech-Gysin and Cech restriction morphisms.

\ms\nin
{\bf 1.5.~Proposition.} {\it Let $H^j(U_{\infty})_{\la}$ denote the
$\la$-eigenspace of $H^j(U_{\infty},\C)$, and similarly for
$H_c^j(U_{\infty})_{\la}$, etc. Let $\nu_{f_U,\la}^j$,
$\nu_{c,f_U,\la}^j$ be as in the introduction.
Then we have for $j\in[0,n]$
$$\hskip-20pt\nu_{f_U,\la}^j=\dim\Gr^W_0H^j(U_{\infty})_{\la}
=\dim\Gr^W_0H^j(X_{\infty})_{\la},
\leqno(1.5.1)$$
\vskip-15pt
$$\q\aligned\nu_{c,f_U,\la}^j&=\dim\Gr^W_{2j}H_c^j(U_{\infty})_
{\la}\\ &=\dim\Ker\bl(\Gr^W_0H^j(X_{\infty})_{\la}\to\mopl_k\,
\Gr^W_0H^j(D_{k,\infty})_{\la}\br),\endaligned
\leqno(1.5.2)$$
where the last morphisms are induced by the restriction morphisms
for $D_k\into X$.}

\ms\nin
{\it Proof.} This follows from the spectral sequences in (1.4.1).
Let $L$ denote the increasing filtration on $H^j(U_{\infty})$,
$H_c^j(U_{\infty})$ associated with the spectral sequences and shifted
by $j$ so that $L$ is the limit of the weight filtration on $H^j(U_t)$,
$H_c^j(U_t)$ for $t\in\De^*$, and
$${}_{\infty}E_2^{-i,j+i}=\Gr^L_{j+i}H^j(U_{\infty}),\q
{}_{\infty}\!{}^cE_2^{i,j-i}=\Gr^L_{j-i}H^j(U_{\infty}).
\leqno(1.5.3)$$
The $E_1$-differentials are induced by the Gysin or restriction
morphisms, and are limits of morphisms of pure Hodge structures of
the same weight. Hence they preserve the center of the symmetry of
the action of $N$, which coincides with the weight of the pure Hodge
structure before taking the limit. Set $d_J:=n-|J|=\dim D_J$.
It is well-known that
$$wt\bl(H^j(D_{J,\infty})\br)\subset\begin{cases}[0,2j]&\h{if}\,\,\,
j\in[0,d_J],\\ [2j-2d_J,2d_J]&\h{if}\,\,\,j\in[d_J,2d_J],\end{cases}
\leqno(1.5.4)$$
where the left-hand side is the set of weights of $H^j(D_{J,\infty})$.
This can be shown by using the invariance of the dimension of the
graded pieces of the Hodge filtration by passing to the limit mixed
Hodge filtration $F$ since the latter together its conjugate Hodge
filtration $\overline{F}$ determines the limit mixed Hodge numbers,
see [De1].

We first show (1.5.2). Using (1.4.1), (1.5.3) and (1.5.4), we get the
fist equality of (1.5.2), since
$$\nu_{c,f_U,\la}^j\le\dim\Gr^W_{2j}H^j_c(U_{\infty})_{\la}=
\dim\Gr^W_{2j}\Gr^L_jH^j_c(U_{\infty})_{\la}\le\nu_{c,f_U,\la}^j.$$
Here the first inequality follows from
$$wt\bl(H^j_c(U_{\infty})\br)\subset[0,2j],$$
the middle equality follows from
$$\Gr^W_{2j}\Gr^L_iH^j_c(U_{\infty})=0\q\h{for}\,\,\,i\ne j,$$
and the last inequality follows from the fact that the
$E_1$-differential preserves the center of the symmetry of the action
of $N$. Moreover, the $E_1$-differential
${}_{\infty}\!{}^cE_1^{0,j}\to{}_{\infty}\!{}^cE_1^{1,j}$ is given by
the restriction morphism
$$H^j(X_{\infty})\to\mopl_k\,H^j(D_{k,\infty}).$$
So we get also the second equality of (1.5.2).

The argument is similar for (1.5.1), and is simpler since we use
in this case the Gysin morphism
$$\mopl_k\,H^{j-2}(D_{k,\infty})(-1)\to H^j(X_{\infty}),$$
where the image has weights in $[2,2j-2]$ so that it can be neglected
for the calculation of $\nu_{f_U,\la}^j$.
This finishes the proof of Proposition~(1.5).

\ms\nin
{\bf 1.6.~Remark.} If we replace the complex manifold $X$ with a
K\"ahler manifold $X'$ having a bimeromorphic proper morphism
$X'\to X$, then $\nu_{f,\la}$ does not change. Indeed, $H^j(X_t,\Q)$
is a direct factor of $H^j(X'_t,\Q)$ for $t\in\De^*$, and the level
of its complement is strictly less than $\min(j,2\dim X_t-j)$.
Here the level of a mixed Hodge structure $H$ is the difference
between the maximal and minimal integers $p$ with $\Gr^p_FH_{\C}\ne 0$.

\ms\nin
{\bf 1.7.~Proof of Theorem~2.}
We can define the spectral sequence (1.3.1) together with the
decomposition (1.3.2) without assuming $X$ K\"ahler.
We have to show its $E_2$-degeneration together with the symmetry of
the $E_2$-term by the action of $N$ (i.e. the induced filtration on
$H^j(X_{\infty})$ is the monodromy filtration with center $0$).
By hypothesis, there is a proper surjective morphism from a K\"ahler
manifold $X'$ to $X$.
Then, using the decomposition theorem for $X'\to X$ (see [Sa3]), the
above properties are reduced to the K\"ahler case, and then follows
from [S1], Section~4 (or [GuNa]).
So the assertion in the case $D=\emptyset$ follows from (1.3.5) by
setting
$$B_{f,\la}^j:=H^0(Y,\Cc_{Y,\la}^j).
\leqno(1.7.1)$$
The general case is then reduced to the case $D=\emptyset$ by
Proposition~(1.5). This completes the proof of Theorem~2.

\ms\nin
{\bf 1.8.~Proof of Theorem~3.}
The nearby and vanishing cycle functors commute with the direct image
by the proper morphism $f':X'\to\De$. So we get
$$\varphi_t\R f'_*\C_{X'}[n]=(\varphi_{f'}\C_{X'}[n])_0,
\leqno(1.8.1)$$
where the right-hand side is identified with the reduced Milnor
cohomology at $0\in X'$ (which is the only singular point of $f'$).
We have furthermore
$$\varphi_{t,\la}\R f'_*\C_{X'}[n]=\begin{cases}\psi_{t,\la}R^nf'_*\C_{X'}
=\psi_{t,\la}R^nf_*\C_X&\h{if}\,\,\,\la\ne 1,\\
{\rm Im\,can}\oplus{\rm Ker\,var}&\h{if}\,\,\,\la=1,\end{cases}
\leqno(1.8.2)$$
where ${\rm can}:\psi_{t,1}\to\varphi_{t,1}$ and
${\rm var}:\varphi_{t,1}\to\psi_{t,1}(-1)$ are as in [Sa1], Section 5,
and we apply these to $^p\HH^0\R f'_*(\C_{X'}[\dim X'])$
(see [BBD] for $^p\HH^j$).
The assertion for $\la=1$ follows from the decomposition theorem in
loc.~cit.
We have moreover
$${\rm Im\,can}={\rm Im}\,N\subset\psi_{t,1}R^nf'_*\C_{X'}=
\psi_{t,1}R^nf_*\C_X,
\leqno(1.8.3)$$
and the action of $N$ on ${\rm Ker\,var}$ is trivial.
We thus get for any $\la$
$$\nu_{g_0,\la}^n=\nu_{f,\la}^n.
\leqno(1.8.4)$$
(Here it is not necessary to assume that the restriction morphism
induces a surjection from $H^n(X_t,\C)$ to the Milnor cohomology.)

On the other hand, we have
$$\varphi_t\,{}^p\HH^j\R f'_*(\C_{X'}[\dim X'])=0\q\h{if}\,\,\,j\ne 0,
\leqno(1.8.5)$$
since $f'$ has only isolated singularities and the vanishing cycle
functor commutes with the direct image by proper morphisms.
This implies that the local systems
$$R^jf'_*\C_{X'}|_{\De^*}=R^jf_*\C_X|_{\De^*}$$
are constant for $j\ne n$, and hence $\nu_{f,\la}^j=0$ if
$j\in[1,n-1]$ or $j=0$ with $\la\ne 1$, where $\nu_{f,1}^0=1$. 
So the assertion follows from Theorem~2.

\bs\bs
\centerline{\bf 2. Partial generalization to the singular case}
\bs\nin
In this section we partially generalize Theorem~2 to the singular case
in Theorem~(2.2).
\ms\nin
{\bf 2.1.~Singular case.}
Theorem~2 for $\nu^j_{f_U,\la}$ cannot be generalized to the singular
case, see Example~(2.3) below. However, we can generalize the assertion
for $\nu^j_{c,f_U,\la}$ in Theorem~2 to the singular case as follows.
Let $f:X\to\De$ be a projective morphism of a reduced analytic space $X$
to $\De$, and $D$ be a closed reduced analytic subspace of $X$ such that
any irreducible components of $X$ and $D$ are dominant over $\De$.
Set $U:=X\setminus D$ with $f_U:U\to\De$ the morphism induced by $f$.
Let $n:=\dim X-1$.
\sk
Let $\nu_{c,f_U,\la}^j$ and $\nu_{c,f_U,\la}^{2n-j}$ be respectively
the number of Jordan blocks of size $j$ and eigenvalue $\la$ for the
monodromy on $H_c^j(U_t)$ and $H_c^{2n-j}(U_t)$ with $j\le n$.
For the statement of Theorem~(2.2) below for
$\nu_{c,f_U,\la}^{2n-j}\,\,(j\le n)$, it is enough to take a
resolution of singularities $\pi_{(0)}:X_{(0)}\to X$ with $\pi_{(0)}$
projective.
For $\nu^j_{c,f_U,\la}\,\,(j\le n)$, however, the
preparation for Theorem~(2.2) is more complicated. We have to
construct
complex manifolds $X_{(0)}$, $X_{(1)}$, $D_{(0)}$ together with
projective morphisms $\pi_{(k)}:X_{(k)}\to X\,\,(k=1,2)$,
$\pi'_{(0)}:D_{(0)}\to D$ and an analytic cycle $\ga_X$ on
$X_{(1)}\times X_{(0)}$ which is a $\Z$-linear combination of graphs of
morphisms from connected components of $X_{(1)}$ to $X_{(0)}$ over $X$
(where there may be many morphisms defined on one connected component).
They have to satisfy the following conditions:
\sk\nin
(i) The composition $f_{(k)}:=f\ssc\pi_{(k)}:X_{(k)}\to\De$ is flat,
(i.e. any connected component is dominant over $\De$), its restriction
over $\De^*$ is smooth, and $f_{(k)}^{-1}(0)$ is a divisor with simple
normal crossings on $X_{(k)}\,\,(k=1.2)$.
\sk\nin
(ii) We have $\Ga_{\pi_{(0)}}\ssc\ga_X=0$ as a cycle on $X_{(1)}\times X$
(without any equivalence relation) where $\Ga_{\pi_{(0)}}$ is the graph
of $\pi_{(0)}$, and the composition of correspondences
$\Ga_{\pi_{(0)}}\ssc\ga_X$ is defined in this case by using the
composition of morphisms.
\sk\nin
(iii) Setting $X_t:=f^{-1}(t)$, $X_{(k),t}:=f_{(k)}^{-1}(t)$, we have
the following exact sequence for any $j\in\Z$ and $t\in\De^*$:
$$0\to\Gr^W_jH^j(X_t,\Q)\buildrel{\pi_{(0)}^*}\over\longrightarrow
H^j(X_{(0),t},\Q)\buildrel{\ga_X^*}\over\longrightarrow
H^j(X_{(1),t},\Q),
\leqno(2.1.1)$$
where $W$ is the weight filtration of the canonical mixed Hodge
structure on $H^j(X_t,\Q)$, and $\ga_X^*$ is defined by using the
pull-backs by the morphisms in the definition of $\ga_X$.
\sk\nin
(iv) The above condition~(i) for $k=0$ with $X$ replaced by $D$ is
satisfied, where we denote the restriction of $f$ to $D$ by $h$,
and the morphism $D_{(0)}\to\De$ by $h_{(0)}$. Moreover $\pi'_{(0)}$ is
surjective and there is a morphism $\rho_{(0)}:D_{(0)}\to X_{(0)}$
giving a commutative diagram
$$\begin{matrix}&&D_{(0)}&\buildrel{\pi'_{(0)}}\over\longrightarrow&D\\
&&\,\,\da{\scriptstyle\rho_{(0)}}&&
\,\,\da{\scriptstyle i}\\
X_{(1)}&\buildrel{\ga_X}\over\longrightarrow&X_{(0)}&
\buildrel{\pi_{(0)}}\over\longrightarrow&X\end{matrix}
\leqno(2.1.2)$$
(Here $X_{(1)}$ is noted also since this will be useful for
Theorem~(2.2) below.)
\ms
This can be done for instance by using an argument similar to [GNPP]
together with resolution of singularities. 
If $X$, $D$ are defined algebraically (i.e. if they are base changes
of algebraic varieties over a curve $C$ by an open inclusion
$\De\into C^{\rm an}$), then the above assumptions are satisfied by
using simplicial resolutions [De3] or cubic resolutions [GNPP].
\sk
Let $B_{f_{(k)},\la}^{\ssb}$ be as in Theorem~2 applied to
$f_{(k)}:X_{(k)}\to\De\,\,(k=1,2)$, and similarly for
$B_{h_{(0)},\la}^{\ssb}$. For any morphism $g$ of a connected component
of $X_{(1)}$ to $X_{(0)}$, we have a morphism of complexes
$$g^*:B_{f_{(0)},\la}^{\ssb}\to B_{f_{(1)},\la}^{\ssb},$$
by choosing an irreducible component of $f_{(0)}^{-1}(0)$ containing the
image of each irreducible component of $f_{(1)}^{-1}(0)$ by $g$.
This induces a morphism of complexes
$$\ga_X^*:B_{f_{(0)},\la}^{\ssb}\to B_{f_{(1)},\la}^{\ssb},$$
and similarly for
$$\rho_{(0)}^*:B_{f_{(0)},\la}^{\ssb}\to B_{h_{(0)},\la}^{\ssb}.$$
\ms\nin
{\bf 2.2.~Theorem.} {\it With the above notation and assumptions, we 
have for $j\in[0,n]$}
$$\nu_{c,f_U,\la}^j=\dim\Ker\bl((\ga_X^*,\rho_{(0)}^*):
H^jB_{f_{(0)},\la}^{\ssb}\to H^jB_{f_{(1)},\la}^{\ssb}\oplus
H^jB_{h_{(0)},\la}^{\ssb}\br),
\leqno(2.2.1)$$
\vskip-15pt
$$\nu_{c,f_U,\la}^{2n-j}=\dim H^jB_{f_{(0)},\la}^{\ssb}.
\leqno(2.2.2)$$
\ms\nin
{\it Proof.}
We first consider $\nu_{c,f_U,\la}^j\,\,(j\le n)$, and prove (2.2.1).
We have a long exact sequence of mixed Hodge structures for $t\in\De^*$
$$\to H^{j-1}(D_t,\Q)\to H_c^j(U_t,\Q)\to H^j(X_t,\Q)
\buildrel{i^*}\over\to H^j(D_t,\Q)\to.
\leqno(2.2.3)$$
Since $H^{j-1}(D_t)$ has weights at most $j-1$, this induces an
isomorphism
$$\Gr^W_jH_c^j(U_t,\Q)=
\Ker\bl(i^*:\Gr^W_jH^j(X_t,\Q)\to\Gr^W_jH^j(D_t,\Q)\br).$$
Combining this with (2.1.1) and using (2.1.2), we see that
$\Gr^W_jH_c^j(U_t,\Q)$ is isomorphic to the kernel of
$$(\ga_X^*,\rho_{(0)}^*):H^j(X_{(0),t},\Q)\to H^j(X_{(1),t},\Q)
\oplus H^j(D_{(0),t},\Q).$$
By [De3] (and using (2.2.3)), we have
$$\Gr_F^pW_{j-1}H_c^j(U_t,\Q)=0\q\h{for}\,\,p\notin[0,j-1].$$
So the assertion (2.2.1) follows from the same argument as in the proof
of Proposition~(1.5).
\sk
We now consider $\nu_{2n-c,f_U,\la}^j\,\,(j\le n)$.
For the proof of (2.2.2), note first that we may replace $X_{(0)}$
by any resolution of singularities of $X$ by Remark~(1.6).
Here we can neglect any complex manifold $Y$ of pure dimension $m<n$,
since we use duality and the dual of $\Q_Y$ is
$$\Q_Y(m)[2m]=\bl(\Q_Y(m)[2n]\br)[2r]\q\h{with}\,\,\,r:=m-n<0.$$
(Without using duality, it is related to the fact that the level of
$H^j(Y,\Q)$ is strictly less than $j$ if $j>\dim Y$.)
We can construct also $X_{(1)}$ and $D_{(0)}$ as in the above case
by using an argument as in [GNPP] so that we may assume moreover
$$\dim X_{(1)}<n,\q\dim D_{(0)}<n.$$
Using the dual argument of the proof of (2.2.1), we get only the
Cech-Gysin morphisms. So (2.2.2) follows from the same argument as in
the proof of Proposition~(1.5).
This finishes the proof of Theorem~(2.2).
\ms
The following example shows that Theorem~2 for $\nu_{f_U,\la}$ cannot
be generalized to the singular case.

\ms\nin
{\bf 2.3.~Example.} Let $Z'=\PP^1$ with $\Sigma:=\{0,\infty\}$.
Let $\sigma:Z'\to Z$ be a morphism inducing an isomorphism outside
$\Sigma$, and such that $\sigma(\Sigma)$ is one point.
Let $\iota':\De\into\PP^1$ be the natural inclusion of an open
disk $\De$ of radius $<1$. This induces an inclusion
$\iota:\De\into Z$, and $1\in\PP^1\setminus\iota'(\De)$ is
identified with a point of $ Z\setminus\iota(\De)$ which is also
denoted by 1. Set $X:=Z\times\De$ with $f:X\to\De$ the second projection.
Let $D\subset X$ be the union of the graph of $\iota$ and
$\{1\}\times\De$.
Set $U:=X\setminus D$. Then, for $t\in\Delta^*$, we have isomorphisms
$$H^1(U_t)=H^1(Z\setminus\{1,t\},\sigma(\Sigma))=
H^1(\PP^1\setminus\{1,t\},\Sigma),$$
where the cohomology is with $\Q$-coefficients, and $\iota'(t)$,
$\iota(t)$ are denoted by $t$ to simplify the notation.
We have a long exact sequence
$$H^0(\PP^1\setminus\{1,t\})\to H^0(\Sigma)\to
H^1(\PP^1\setminus\{1,t\},\Sigma)\to
H^1(\PP^1\setminus\{1,t\})\to 0,$$
inducing a short exact sequence
$$0\to\Q\to H^1(\PP^1\setminus\{1,t\},\Sigma)\to\Q(-1)\to 0.$$
We see that the monodromy around the origin in $\De$ is nontrivial as
follows. There is a relative cycle class $\ga$ in 
$H_1(\PP^1\setminus\{1,t\},\Sigma)$ represented by a path between
$0$ and $\infty$ which is slightly below the real positive half line.
Let $t_0\in\De^*$ be a sufficiently small real positive number.
Take a loop $\al\in\pi_1(\De^*,t_0)$ going around the origin of $\De$
counterclockwise. Deform the relative cycle $\ga$ continuously when
$t\in\De^*$ moves along $\al$. Then the relative cycle $\ga$
becomes slightly above the real positive half line. Thus the action of
the monodromy $T$ on the relative cycle $\ga$ is given by
$$T\ga=\ga+\eta,$$
where $\eta$ is a small circle around $t_0$.
This implies the non-vanishing of
$$N:H^1(U_{\infty})\to H^1(U_{\infty})(-1).$$
However, we have
$$\Gr^L_kH^1(U_{\infty})=\Gr^W_kH^1(U_{\infty})=\begin{cases}\Q &
\h{if}\,\,\,k=0,\\ \Q(-1) &\h{if}\,\,\,k=2,\\ \,0 &\h{if}\,\,\,
k\ne 0,2,\end{cases}$$
where $L$ is induced by the weight filtration $W$ on $H^1(U_t)$ for
$t\in\De^*$ as in the proof of Proposition~(1.5).
Thus Theorem~2 for $\nu_{f_U,\la}$ is false in the singular case.

\bs\bs
\centerline{\bf 3. Global triviality of certain nearby cycle local
systems}
\bs\nin
In this section we provide a method to show Theorem~1 in Section~4,
and prove Theorem~4.
\ms\nin
{\bf 3.1.~Global factorization of functions.}
Let $f$ be a holomorphic function on a complex manifold.
Assume $Y:=f^{-1}(0)$ is a divisor with simple normal crossings.
Set $X^*:=X\setminus Y$ with the inclusion $j:X^*\into X$.
For a locally closed analytic subset $Z$ of $Y$ with the inclusion
$i_Z:Z\into X$, set
$$\M^*_Z:=i_Z^{-1}j_*^{\rm mer}\OO_{X^*}^*,$$
where $j_*^{\rm mer}$ denotes the meromorphic extension over $X$.
If $Y$ is locally the union of $\{x_i=0\}$ for $0\le i<r$, where
$x_0,\dots,x_n$ are local coordinates of $X$ around $x\in Y$, then
$$\M^*_{Z,x}=\bl\{u\,\mprod_{0\le i<r}\,x_i^{a_i}\,\big|\,
u\in\OO_{X,x}^*,\,a_i\in\Z\br\}.$$
For an integer $m'\ge 2$, let $\M^{*m'}_Z$ be the image of the $m'$-th
power endomorphism of $\M^*_Z$. We have a short exact sequence of
sheaves of multiplicative groups over $Z$
$$1\to\mu_{Z,m'}\to\M^*_Z\buildrel{m'}\over\longrightarrow\M^{*m'}_Z
\to1,$$
where $\mu_{Z,m'}$ is the constant sheaf on $Y$ with stalks $\mu_{m'}$
(the multiplicative group consisting of the roots of unity of order
$m'$ in $\C^*$).
We have the associated long exact sequence
$$1\to\mu_{m'}\to\Ga(Z,\M^*_Z)\buildrel{m'}\over\longrightarrow
\Ga(Z,\M^{*m'}_Z)\buildrel{c_{m'}}\over\longrightarrow
H^1(Z,\mu_{Z,m'}),$$
where the last morphism $c_{m'}$ gives the cohomology class of
$u\in\Ga(Z,\M^{*m'}_Z)$. This is the same as the cohomology class of
the finite unramified covering of $Z$ defined by $m'{}^{-1}(u)$
which is a principal $\mu_{m'}$-bundle.
(Indeed, consider the Cech cocycle associated to local pull-backs of
$u$ by $m'$ for an sufficiently fine open covering of $Z$.)
Anyway, we have a primitive $m'$-th root of $u$ globally over $Z$
if and only if $c_{m'}(u)=0$.

Assume the restriction of $f$ to a sufficiently small neighborhood
of $Z$ defines an element $u_f$ of $\Ga(Z,\M^{*m'})$, i.e. there is a
solution of $\xi^{m'}=f$ with $\xi\in\OO_{X,x}$ for any $x\in Z$.
Then we have a global solution of $\xi^{m'}=f$ on a sufficiently small
open neighborhood of $Z$ if and only if $c_{m'}(u_f)=0$.

\ms\nin
{\bf 3.2.~Globally factorized case.}
With the notation of (3.1), assume there is a global solution
$\xi^{m'}=f$ on $X$ where $f:X\to\De$ is not necessarily proper.
We have a factorization
$$f:X\buildrel{\xi}\over\longrightarrow\Det'
\buildrel{\pi_{m'}}\over\longrightarrow\De,$$
where $\pi_{m'}$ is a totally ramified covering of degree $m'$.
Let $\Xt'$ be the normalization of the base change of $f:X\to\De$ by
$\pi_{m'}:\Det'\to\De$. Let $\pi':\Xt'\to X$ be the canonical morphism.
Set $\Yt':=\pi'{}^{-1}(Y)$ with $\pi'_0:\Yt'\to Y$ the canonical
morphism. Note that this is a trivial covering space, i.e.
$\Yt'$ is a disjoint union of $m'$ copies of $Y$.

Let $V_{m'}$ be a complex vector space endowed with a basis
$(e_0,\dots,e_{m'-1})$ and an action of $T$ defined by $Te_i=e_{i+1}$
for $i=0,\dots,m'-1\,\,\,\h{mod}\,\,\,m'$.
Let $V_{m',Y}$ denote the constant sheaf with stalks $V_{m'}$.
Then, choosing a section of $\pi'_0$, we have canonical isomorphisms
$$\mopl_{\la^{m'}=1}\,\HH^0\psi_{f,\la}\C_X=(\pi'_0)_*\C_{\Yt'}
=V_{m',Y},
\leqno(3.2.1)$$
in a compatible way with the action of $T$,
where $T$ on the middle term is given by the action of an appropriate
generator of the covering transformation group of $\pi'_0$.
Indeed, the first isomorphism is shown by using the Milnor fiber at
each point. The second isomorphism follows from the triviality of
the covering $\pi'_0:\Yt'\to Y$ by choosing a section of $\pi'_0$.

For $I\subset J(\la)$, set
$$L'_{\la,I}:=L_{\la,I}|_{U'_{\la,I}}\q\h{with}\q
U'_{\la,I}:=U_{\la,I}\cap Y^{(\la)}.$$
Then, using the projection from $V_{m'}$ to $\C\,e_0\subset V_{m'}$,
(3.2.1) induces canonical isomorphisms
$$L'_{\la,I}=\C_{U'_{\la,I}},
\leqno(3.2.2)$$
in a compatible way with the restriction morphisms (1.2.8).

\ms\nin
{\bf 3.3.~Proposition.} {\it For an integer $m'\ge 2$, let $Z$ be a
closed subvariety of $Y^{(\la')}$ with $\la':=\exp(2\pi i/m')$ in the
notation of Theorem~$4$. Let $\pi_Z:\Zt\to Z$ be the unramified
covering of degree $m'$ defined by local solutions of $\xi^{m'}=f$ as
in $(3.1)$. Then, with the notation of $(3.2.2)$, we have a canonical
isomorphism
$$(\pi_Z)_*\C_{\Zt}\simto
\mopl_{\la^{m'}=1}\,\bl(\HH^0\psi_{f,\la}\C_X\br)\big|_Z
\leqno(3.3.1)$$
in a compatible way with the action of $T$ where $T$ on the
left-hand side is defined by the action of an appropriate generator
of the covering transformation group of $\pi_Z$.}

\ms\nin
{\it Proof.} Let $\Xt'$ be the normalization of the base change of
$f:X\to\De$ by the $m'$-fold ramified covering $\pi_{m'}:\Det'\to\De$.
Restricting over a sufficiently small open neighborhood of each $z$ of
$Z$ in $X$, this coincides with the construction in (3.2).
Note that the restriction of $\Xt'\to X$ over $Z$ is identified with
$\pi_Z:\Zt\to Z$. Let $\ft':\Xt'\to\Det'$ be the natural morphism.
We have natural isomorphism and inclusion
$$(\pi_Z)_*\C_{\Zt}\simto
(\pi_Z)_*\bl(\HH^0\psi_{\ft',1}\C_{\Xt'}\big|_{\Zt}\br)\into
\HH^0\psi_f\C_X\big|_Z$$
compatible with the action of $T$, where $T$ on the fist and second
terms is induced by the action of an appropriate generator of the
covering transformation group of $\pi_Z$. So the assertion follows
from the local calculation in (3.2.1) (by counting the dimension).
This finishes the proof of Proposition~(3.3).

\ms\nin
{\bf 3.4.~Corollary.} {\it With the above notation and assumption,
$\pi_Z:\Zt\to Z$ is trivial, if and only if the $\HH^0\psi_{f,\la}\C_X$
are trivial local systems for any $\la$ with $\la^{m'}=1$.}

\ms\nin
{\it Proof.} This follows from Proposition~(3.3) by applying the
global section functor to (3.3.1).

\ms\nin
{\bf 3.5.~Proof of Theorem~4.}
We apply (3.1) to the case $Z=Y^{(\la)}$ and $m'=m_{\la}$.
If $H^1(Y^{(\la)},\mu_{m_{\la}})=0$, then we have a global
solution of $\xi_{\la}^{m_{\la}}=f$ on a sufficiently small open
neighborhood $X^{(\la)}$ of $Y^{(\la)}$.
So the assertion (ii) follows from (3.2).
The argument is similar for the remaining assertions.
This finishes the proof of Theorem~4.

\ms\nin
{\bf 3.6.~Proposition.} {\it With the notation of Theorem~$4$, assume
there is a subset $Z^{(\la)}$ of $Y^{(\la)}$ which is homotopy
equivalent to a dense Zariski-open subset $U^{(\la)}$ of $Y^{(\la)}$,
and moreover there is a holomorphic function $g_{\la}$ on a
sufficiently small open neighborhood of $Z^{(\la)}$ in $X$ satisfying
$g_{\la}^{m_{\la}}=f$ on this neighborhood.
Then $B_{f,\la}^{\ssb}=C_{f,\la}^{\ssb}$.
If the above condition holds by replacing $X$, $Y^{(\la)}$ and $f$
respectively with $D_k$, $Y_k^{(\la)}$ and $f_k=f|_{D_k}$ for any $k$,
then $B_{f_k,\la}^{\ssb}=C_{f_k,\la}^{\ssb}$.}

\ms\nin
{\it Proof.}
By the same argument as in the proof of Theorem~4, it is sufficient
show that we have a global solution of $\xi^{m_{\la}}=f$ on a
sufficiently small open neighborhood $X^{(\la)}$ of $Y^{(\la)}$.
Here we may replace $Y^{(\la)}$ with the dense Zariski-open subset
$U^{(\la)}$.
Indeed, local solutions of $\xi^{m_{\la}}=f$ form a finite unramified
covering as in (3.1), and it is trivial over $Y^{(\la)}$ if its
restriction over any dense Zariski-open subset is trivial.
Moreover, the triviality of the covering is determined by its
cohomology class in the first cohomology with coefficients in
$\mu_{m_{\la}}$, see (3.1).
This triviality can be seen by restricting to the subspace $Z^{(\la)}$
which is homotopy equivalent to $U^{(\la)}$ by the hypothesis of
Theorem~4. So the assertion follows.

\ms\nin
{\bf 3.7.~Proposition.} {\it With the notation and the assumption
of Theorem~$3$, assume $n=2$ and the embedded resolution is obtained by
iterating blowing-ups with point or $\PP^1$-centers.
Then $B_{f,\la}^j=C_{f,\la}^j$ for any $j$, and hence Theorem~$3$ holds
with $B_{f,\la}^{\ssb}$ replaced by $C_{f,\la}^{\ssb}$.}

\ms\nin
{\it Proof.} Since projective spaces $\PP^k\,(k=1,2)$ and $\PP^1$-bundles
over $\PP^1$ are simply connected, and simple connectedness does not
change by point-center blow-ups, the assertion follows from Theorem~4(i).
\ms
The following is closely related with results in [Ar1], [MM]
where similar constructions are used.

\ms\nin
{\bf 3.8.~Proposition.} {\it With the notation and the assumption
of Theorem~$3$, assume $n=2$ and $g_0$ defines a super-isolated
singularity {\rm [Lu]} or more generally, a Yomdin singularity {\rm [Yo]}.
Then $B_{f,\la}^j=C_{f,\la}^j$ for any $j$, and hence Theorem~$3$ holds
with $B_{f,\la}^{\ssb}$ replaced by $C_{f,\la}^{\ssb}$.}

\ms\nin
{\it Proof.} We have the expansion $g_0=\sum_{j\ge d}g_{0,j}$ with
$g_{0.j}$ a homogeneous polynomial of degree $j$, and $g_{0,d}\ne 0$. Set
$$Z:=g_{0,d}^{-1}(0)\subset\PP^2.$$
Then the condition that $g_0^{-1}(0)$ is a Yomdin singularity [Yo]
means that $Z$ has only isolated singularities, $g_{0,j}=0$ for $d<j<k$,
and $g_{0,d+k}^{-1}(0)\cap{\rm Sing}\,Z=\emptyset$, see [ALM].
It is a super-isolated singularity [Lu] if $k=1$.
We show that the embedded resolution can be obtained by repeating
blowing-ups with point or $\PP^1$-centers.

We first take the blow-up $\sigma_1:X_1\to X_0=X'$ at $0\in X'$.
Its exceptional divisor $E_0$ is $\PP^2$, and the intersection of $E_0$
with the proper transform of $g_0^{-1}(0)$ is identified with
$Z\subset\PP^1$.
Moreover, the total transform of $g_0^{-1}(0)$ around a singular point of
$Z$ can be defined locally by an equation of the form
$$w^d(h(u,v)+w^k)=0,
\leqno(3.8.1)$$
where $(u,v,w)$ is a local coordinate system such that the exceptional
divisor $E_0=\PP^2$ is locally defined by $w=0$, and $Z\subset\PP^2$ is
defined by $h(u,v)=0$. Here the restrictions of $x,y$ to $\PP^2$ are
identified with local coordinates of $\PP^2$.
Indeed, take a coordinate system $(x,y,z)$ of $\C^3$.
Set $h_j:=g_{0,j}/z^j$. This is viewed as a function on the complement of
$\{z=0\}\subset E_0=\PP^2$. Then the pull-back of $g_0$ to the complement
of the proper transform of $\{z=0\}\subset\C^3$ is expressed as
$$z^d\bl(h_d+z^k(h_{d+k}+zh_{d+k+1}+\cdots)\br),$$
where $z$ denotes also the pull-back of $z$ which locally defines the
exceptional divisor $E_0$. So (3.8.1) follows by setting locally
$$h:=h_d,\q w:=z(h_{d+k}+zh_{d+k+1}+\cdots)^{1/k}.$$

Repeating point-center blow-ups at singular points of the total transform
of $Z$ in the proper transform of $E_0=\PP^2$, we then get a morphism
$$\sigma_2:X_2\to X_1,$$
such that the intersection of the total transform of $Z$ with the
proper transform $\E_0$ of $E_0$ is a divisor with simple normal
crossings on $\E_0$.
(Here we use the fact that the restriction of a point-center blow-up
to the proper transform of a smooth divisor is a point-center blow-up.)
We may moreover assume that any two irreducible components of the
proper transform of $Z$ do not intersect each other
(taking a point-center blow-up at the intersection point if necessary).

Applying a point-center blow-up to (3.8.1), the local coordinate system
$(u,v,w)$ is substituted by $(u,uv,uw)$ or $(uv,v,vw)$ near the proper
transform of $E_0$.
Repeating this, the total transform of $g_0^{-1}(0)$ by
$\sigma_1\ssc\sigma_2$ is locally defined by
$$u^iv^jw^l(u^av^b+w^c)=0\q\h{with}\q i,j\ge 0,\,\,l,a,b,c>0,
\leqno(3.8.2)$$
using a local coordinate system $(u,v,w)$, where $l=d$, $c=k$.
We have $a=1$ if $i=0$, and $b=1$ if $j=0$.
Note that the non-normal crossing points of (3.8.2) are contained in the
union of $\{u=w=0\}$ and $\{v=w=0\}$.

By the above construction, the non-normal crossing points of the total
transform of $g_0^{-1}(0)$ consist of a union of smooth rational curves.
In order to apply Proposition~(3.7), it is then sufficient to show
that (3.8.2) is essentially {\it stable} by blowing-ups along the origin
or along the coordinate axes.
Indeed, it is known that Hironaka's resolution can be obtained by
repeating blow-ups with smooth centers contained in the set of non-normal
crossing points, and the {\it new components} of the set of non-normal
crossing points which are obtained by a blow-up of the divisor defined by
the equation of the form (3.8.2) are also {\it rational curves}.
Here ``essentially" means that we allow $a,b,c\ge 0$ together with a
{\it permutation of variables} and that we may get an equation which is
not of the form (3.8.2) if the equation defines a divisor with normal
crossings as explained below.
(It may be possible to give a more explicit algorithm by induction on the
maximum of $a,b,c$, although this seems more complicated than one might
imagine. Indeed, a resolution of singularities is {\it global} on $X_2$
and a {\it local description} using the Euclidean algorithm at each point
of $X_2$ is not enough. Here permutations of variables make the argument
rather complicated.)

In case of a point-center blow-up, $(u,v,w)$ is substituted in (3.8.2) by
$$(u,uv,uw)\,\,\,\,\h{or}\,\,\,\,(uv,v,vw)\,\,\,\,\h{or}\,\,\,\,(uw,vw,w).$$
In case of the blow-up along $\{u=w=0\}$, $(u,v,w)$ is substituted in
(3.8.2) by
$$(u,v,uw)\,\,\,\,\h{or}\,\,\,\,(uw,v,w).$$
and similarly for $\{u=w=0\}$ with $u$ replaced
by $v$. By these substitutions, (3.8.2) is essentially stable except for
the case we get a local equation of the form
$$u^iv^jw^l(u^av^bw^c+1)=0.
\leqno(3.8.3)$$
However, this defines a divisor with normal crossings, and we do not have
to consider it. So the assertion follows from Proposition~(3.7).
More precisely, under a substitution by $(u,uv,uw)$ or $(u,v,uw)$
for instance, only $a$ changes and $b$, $c$ do not change in (3.8.2)
if we allow $a$ negative, and we need a permutation of variables if we
want to get $a,b,c\ge 0$.
We may have (3.8.3) under a substitution by $(uw,vw,w)$ or $(uw,v,w)$ or
$(u,vw,w)$.
Then, repeating the blow-ups consisting of Hironaka's resolution
(or using the Euclidean algorithm essentially), we will reach local
equations of the form
$$u^iv^jw^l(v^c+w^c)=0\q\h{or}\q v^jw^l(uv^c+w^c)=0\q\h{with}\q c\ge 1,
\leqno(3.8.4)$$
just before getting a divisor with normal crossings by blowing-up along
$\{v=w=0\}$. Here the obtained equation depends on whether we started
from (3.8.2) with $i,j\ge 1$, $ab\ne 0$ or not.
In the latter case, if we start from (3.8.2) with $i=0$, $a=1$, then
we can apply the Euclidean algorithm to $b,c$ in (3.8.1) since
we get a divisor with normal crossings by an equation of the form
$$v^jw^l(u+v^bw^c)=1.$$
This finishes the proof of Proposition~(3.8).

\ms\nin
{\bf 3.9.~Remark.} In case of super-isolated singularities, or more
generally, Yomdin singularities with $n\ge 2$, formulas are known for
the Milnor number, the characteristic polynomial of the Milnor monodromy,
and also for the spectrum, see [ALM], [LM], [Si], [Stv], [Yo].

In fact, Steenbrink ([St3], Th.~6.1) proved a formula for the spectrum
of a homogeneous polynomial $f$ with one-dimensional singular locus,
which can be expressed for instance (using the normalization as in [Sa4])
as follows:
$$\Sp(f,0)=\biggl(\frac{t-t^{1/d}}{t^{1/d}-1}\biggr)^{n+1}-
\sum_{i,j}t^{\al'_{i,j}}\frac{t-1}{t^{1/d}-1},
\leqno(3.9.1)$$
where $\al'_{i,j}:=(\lfloor\al_{i,j}d\rfloor+1)/d$ with
$\lfloor\al\rfloor:=\max\{p\in\Z\mid p\le\al\}$, and the $\al_{i,j}$ are
the exponents, i.e. the spectral numbers counted with multiplicities at
each singular point $y_i$ of $f^{-1}(0)\subset\PP^n$.
Note that (3.9.1) is quite useful for calculations of the spectrum in
this case; for instance, the formula in [BS], Th.~3 for the spectrum of
reduced hyperplane arrangements in $\C^3$ follows from it.

We may view (3.9.1) as a special case (with $k=0$) of Steenbrink's
conjecture in [St3], which was proved there in case $f$ is homogeneous
and the isolated singularities are of Brieskorn type (and in [Sa4] in
general). The latter can be expressed in this case as follows:
$$\Sp(f+h^{d+k},0)-\Sp(f,0)=
\sum_{i,j}t^{\al''_{i,j}(k)}\frac{t-1}{t^{1/d+k}-1}\q\q(k\ge 0),
\leqno(3.9.2)$$
where $\al''_{i,j}(k):=(k\al_{i,j}+\lfloor\al_{i,j}d\rfloor+1)/(d+k)$, and
$h$ is a sufficiently general linear function, see also [ALM], Th.~1.4.
Indeed, in the homogeneous polynomial case, there is a well-known
relation between the Milnor monodromy and the local system monodromy
along $\C^*\subset{\rm Sing}\,f^{-1}(0)$ so that $\beta_{i,j}$ in
[Sa4], (0.1) satisfies the relation
$$\al_{i,j}d+\beta_{i,j}\in\Z.
\leqno(3.9.3)$$
Combining this with the condition $\beta_{i,j}\in(0,1]$, we get
$$\al_{i,j}d+\beta_{i,j}=\lfloor\al_{i,j}d\rfloor+1,\q\h{and}\q
\al''_{i,j}(k)=((d+k)\al_{i,j}+\beta_{i,j})/(d+k).
\leqno(3.9.4)$$
The lower bound of $k$ in (3.9.2) is $0$, since the number $R$
in [Sa4], Th.~2.5 is $d$ in this case.
(This can be shown by using the natural $\C^*$-action.)

Note that (3.9.1-2) imply a formula for the spectrum of Yomdin
singularities as in [ALM], Th.~1.4 (using the constancy of the spectrum
by $\mu$-constant deformations).
We can verify that the normalization of the formulas (3.9.1-2) is correct,
for instance, in a simple case where $f:=xyz$ (i.e. of type
$T_{\infty,\infty,\infty}$, see [St3]) with $n=2$, $d=3$, and
$f':=f+x^p+y^p+z^p$ (i.e. of type $T_{p,p,p}$) for $p=k+3>3$.
In this case, we have $\al_{i,1}=\beta_{i,1}=1$ for $i=1,2,3$, and
$$\Sp(f,0)=t-2t^2,\q \Sp(f',0)=\Sp(f,0)+3\,\msum_{l=1}^p\,t^{1+l/p}.$$
(There is a shift by one between the normalizations of the spectrum in
[St3] and in [Sa4].)

Since Steenbrink's conjecture is generalized to the case of spectral
pairs [NS], it would imply a certain formula for the number of Jordan
blocks of the Milnor monodromy of Yomdin singularities by using the
monodromical property of the weight filtration [St2].

\ms\nin
{\bf 3.10.~A criterion.}
In the case of Theorem~3, we can determine whether the equality
$B_{f,\la}^j=C_{f,\la}^j$ for $\la\ne 1$ holds in certain cases as
follows.
Here we consider a slightly more general situation where $f:X\to\De$
is obtained by an embedded resolution of the singular fiber
$f'{}^{-1}(0)$ of a morphism of complex manifolds $f':X'\to\De$ where
the singularities of $f'{}^{-1}(0)$ are not necessarily isolated.
We assume the resolution is given by the composition of blow-ups
with connected smooth centers
$$\sigma_i:X_i\to X_{i-1}\q(i=1,\dots,r)$$
where $X_0=X'$ and $X_r=X$.
Let $E_i\subset X_i$ be the exceptional divisor of $\sigma_i$
with $D_i$ its proper transform in $X$.
Let $m_i$ be the multiplicity of $Y$ along $D_i$.
Let $g_i$ be the pull-back of $f'$ to $X_i$.

Fix some $i\in[1,r]$ with $m_i/m_{\la}\in\Z$.
Let $Z$ be a closed subvariety of $D_i\cap Y^{(\la)}$ such that the
canonical morphism $\pi_i:X\to X_i$ induces a morphism of $Z$ to its
image $Z'$ in $X_i$ with {\it connected fibers}.
Assume there is a meromorphic function $h_i$ on a neighborhood $U_{Z'}$
of $Z'\subset X_i$ (in classical topology) satisfying the following three
conditions:

\ms\nin
(i) The zeros of the pull-back of $h_i$ in a sufficiently small open
neighborhood $U_Z$ of $Z$ in $\pi_i^{-1}(U_{Z'})$ are contained in $Y$.

\sk\nin
(ii) The order of zero of $h_i$ along $E_i$ is $m_i/m_{\la}$.

\sk\nin
(iii) The restriction of $g'_i:=g_i/h_i^{m_{\la}}$ to $U_{Z'}\cap E_i$
is a meromorphic function having finite values on dense Zariski-open 
subsets of any intersections of irreducible components of $Z'$.

\ms
Then we have the following (which will be used in (4.3) below).
\ms\nin
{\bf 3.11.~Proposition.} {\it With the above notation and assumption,
there is a global solution of the equation $\xi^{m_{\la}}=f$ on
a sufficiently small neighborhood of $Z$ if and only if there is a
global solution of $\xi'{}^{m_{\la}}=g'_i|_{E_i}$ on $Z'$.}

\ms\nin
{\it Proof.}
Let $g'$ and $h$ respectively denote the pull-back of $g'_i$ and $h_i$
to $U_Z\subset X$. Then
$g'=f/h^{m_{\la}}$, and it is enough to consider the global
solvability of $\xi^{m_{\la}}=g'$.
By hypothesis, the zeros and poles of $g'$ are contained in $Y$,
and it has finite values generically on $U_Z\cap D_i$.
Hence we can take the pull-back of $g'_i$ after restricting it to
$U_{Z'}\cap E_i$.
Then the assertion follows from the hypothesis on the connectivity
of the fibers of the morphism $Z\to Z'$.
This finishes the proof of Proposition~(3.11).

\ms\nin
{\bf 3.12.~Remarks.} (i) In Proposition~(3.11) it is essential to
consider the restriction of $g'_i$ to the intersection with $E_i$,
since $h_i^{-1}(0)$ is not necessarily contained in $g_i^{-1}(0)$ on a
neighborhood of $Z'$ in $X_i$, even though we have the inclusion on a
neighborhood of $Z$ in $X$ after taking the pull-back because of a
blow-up with center contained in the proper transform of
$h_i^{-1}(0)\cap E_i$. This will be used in (4.3).
\ms
(ii) By Proposition~(3.6) for $Z^{(\la)}=Y^{(\la)}$, the global
solvability of the equation $\xi^{m_{\la}}=f$ on a sufficiently small
open neighborhood of $Y^{(\la)}$ implies the equality
$B_{f,\la}^{\ssb}=C_{f,\la}^{\ssb}$.

\bs\bs
\centerline{\bf 4. Examples}
\bs\nin
In this section we give some interesting examples, and prove
Theorem~1 in (4.3).
\ms\nin
{\bf 4.1.~Example.}
Let $E$ be an elliptic curve with the origin $O$.
Let $P$ be a torsion point of $E$ with order $m>1$.
Let $X$ be the blow-up of $E\times E$ along the two points $(O,P)$,
$(P,O)$. Let
$$D_0=E\times\{O\},\q D'_0=\{O\}\times E,\q
D_{\infty}=E\times\{P\},\q D'_{\infty}=\{P\}\times E,$$
and $\Dt_0$, $\Dt'_0$, $\Dt_{\infty}$, $\Dt'_{\infty}$ be their
proper transforms.
Then we have a rational function $f$ on $X$ defining a morphism of
algebraic varieties $f:X\to\PP^1$, and satisfying
$${\rm div}\,f=m\Dt_0+m\Dt'_0-m\Dt_{\infty}-m\Dt'_{\infty}.$$
Indeed, there is a rational function $g$ on $E$ with
${\rm div}\,g=mO-mP$ by Abel's theorem for elliptic curves, and
$f$ is the pull-back of $pr_1^*g\cdot pr_2^*g$ where $pr_1,pr_2$ are
the first and second projections.

However, there is no univalued holomorphic function $g$ with $g^a=f$
for $a>1$ even on a sufficiently small analytic neighborhood of
$f^{-1}(0)$ in $X$ since the general fibers of $f$ are connected.
Indeed, we have finite morphisms $\PP^1\to S\buildrel\rho\over\to\PP^1$
where the first $\PP^1$ is an exceptional divisor of the blow-up, and
$\rho$ is the Stein factorization of $f$.
The composition is given by the restriction of $f$, and is a
ramified covering of degree $m$ which is ramified only at $0$ and
$\infty$.
Then $\rho$ is an isomorphism (i.e. the general fibers of $f$ are
connected), since otherwise there is a rational function $g$ on
$X$ with $g^a=f$ for $a>1$, contradicting the fact that
there is no rational function $g'$ on $E$ with ${\rm div}\,g'=m'O-m'P$
for $0<m'<m$ (by restricting to $E\times\{Q\}$ for a general point
$Q\in E$).

A similar assertion holds by restricting to a neighborhood of
$\Dt_0$ or $\Dt'_0$.
Here we use the first cohomology $H^1(f^{-1}(0),\mu_m)$ as in (3.1).
This gives an example with
$\chi(B_{f,\la}^{\ssb})\ne\chi(C_{f,\la}^{\ssb})$ for
$\la\in\mu_m\setminus\{1\}$.
More precisely, we have for $\la\in\mu_m\setminus\{1\}$
$$B_{f,\la}^0=0,\q C_{f,\la}^0=\C\oplus\C,\q
B_{f,\la}^1=C_{f,\la}^1=\C.$$

In this case, a general fiber $X_t$ is a connected curve of genus
$m+1$ (using for instance the Riemann-Roch theorem on $X$).
Let $H^j(X_{\infty},\Q)$ be the limit mixed Hodge structure, and
$H^j(X_{\infty},\C)_{\la}$ be the $\la$-eigenspace of the monodromy.
Calculating the $E_1$-complex of the weight spectral sequence, we get
$$\Gr^W_kH^j(X_{\infty},\Q)_1=\begin{cases}\Q&\h{if}\,\,\,(j,k)=(0,0),\\
H^1(E,\Q)\oplus H^1(E,\Q)&\h{if}\,\,\,(j,k)=(1,1),\\
\Q(-1)&\h{if}\,\,\,(j,k)=(2,2),\\ \,0&\h{otherwise},\end{cases}$$
and for $\la\in\mu_m\setminus\{1\}$
$$\Gr^W_kH^j(X_{\infty},\C)_{\la}=\begin{cases}\C&\h{if}\,\,\,(j,k)=
(1,0),\\ \C(-1)&\h{if}\,\,\,(j,k)=(1,2),\\ \,0&\h{otherwise}.
\end{cases}$$
In particular, $\nu_{f,1}^1=0$ and $\nu_{f,\la}^1=1$ for any
$\la\in\mu_m\setminus\{1\}$.
(This is the first example with $B^j_{f,\la}\ne C^j_{f,\la}$,
and it was rather surprising.)

\ms\nin
{\bf 4.2.~Example.} Let $C$ be an elliptic curve embedded in
$\PP^2$, and $L_i$ be three lines in $\PP^2$ intersecting $C$ only at
one point $P_i$ with intersection multiplicity 3 for $i=1,2,3$
(i.e. the $P_i$ are inflection points), and such that
$\mcap_{i=1}^3L_i=\emptyset$.
Let $h,h'$ be homogeneous polynomials $h,h'$ of degree 3 defining
$C$ and $\mcup_{i=1}^3\,L_i$ respectively. Using coordinates, we have
$$\aligned&h=x^3+3\al^2x^2y+3\al xy^2+y^3+3(x^2+y^2)z+3(x+y)z^2
+z^3+cxyz,\\ &h'=xyz,\q\h{where $\al^3=1$, and $c\in\C$ is
generic.}\endaligned$$
We assume $\al\ne 1$. This is equivalent to the following condition:
$$\h{The three points $P_1,P_2,P_3$ are {\it not} on the same line in
$\PP^2$.}
\leqno(A)$$
Here we may assume that $P_3$ is the origin $O$ of the elliptic curve.
Then $P_1,P_2$ are torsion points of order 3, and condition~$(A)$ is
equivalent to the condition: $P_1+P_2\ne O$.

Set
$$g':=h^3h':\C^3\to\C.$$
We have an embedded resolution $U'\to\C^3$ of $g'{}^{-1}(0)$ by
blowing-up first the origin, and then repeating the blowing-ups along
the proper transforms of the affine cone of $C\cap L_i$ in $\C^3$
three times for each $i$.
Then the composition $U'\to\C$ of the resolution and $g'$ can be
extended to a projective morphism $f:X\to\C$ such that $X$ is smooth
and $(X\setminus U')\cup f^{-1}(0)$ is a divisor with normal crossings.
However, $U'$ may be different from $U$ in the introduction since
$X\setminus U'$ may contain some vertical divisors.

Let $D'_0\subset U'$ be the proper transform of the exceptional
divisor $\PP^2$ of the first blow-up. Let $D'_i\subset U'$ be the
exceptional divisor of the last blow-up of the successive three
blow-ups along the proper transforms of the affine cone of $C\cap L_i$
for $i=1,2,3$. Let $D'_4\subset U'$ be the proper transform of the
affine cone of $C$.
Let $D_i$ be the closure of $D'_i$ in $X$ for $i=0,\dots,4$ where
$D_0=D'_0$.

Let $m_i$ be the multiplicity of $D_i$.
Then $m_i=12$ for $i=0,\dots,3$, and $m_4=3$.
The multiplicities of the exceptional divisors of the first
and second blow-ups along the proper transforms of the affine cone of
$C\cap L_i$ are respectively 4 and 8, and are not divisible by 3.
So the $D_i$ for $i=0,\dots,4$ are irreducible components of
$f^{-1}(0)$ with multiplicities divisible by 3, and
$D_{\{1,4\}}:=D_0\cap D_4\subset U'$ does not intersect the irreducible
components of $f^{-1}(0)$ other than $D_i\,(i=0,\dots,4)$.
We thus get a unramified covering of degree 3
$$\Dt_{\{1,4\}}\to D_{\{1,4\}},$$
which is non-trivial by condition~$(A)$.
(Indeed, using the coordinates $u=x/z$, $v=y/z$, $w=z$ of the blow-up
at the origin, the pull-back of $g'$ is written as
$\bl(h(u,v,1)w^4\br){}^3uv$. So it is enough to show the non-existence
of a rational function $\xi$ on $C$ satisfying $\xi^3=uv|_C$.
Since ${\rm div}(uv|_C)=3P_1+3P_2-6P_3$, the assertion follows
from the remark after condition~$(A)$.)
We thus get
$$B_{f,\om}^1\ne C_{f,\om}^1\q\h{with}\,\,\,\om=\exp(\pm 2\pi i/3).$$
In this case, the local monodromy is semisimple since $f$ is
homogeneous. In particular, $\nu_{f_U,\la}^j=0$ for $j=1,2$.
This example is needed for the proof of Theorem~1 below.
Note that some related results are obtained in [Ar2], [AC].

\ms\nin
{\bf 4.3.~Proof of Theorem~1.} With the notation of Example~(4.2), set
$$g_0:=h^3h'+h'',$$
where $h''$ is a homogeneous polynomial of degree 16 such that
$h''{}^{-1}(0)\subset\PP^2$ is smooth and transversely intersects
$\mcup_i\,L_i\cup C$ at smooth points.
Let $f$ be a desingularization of a good projective compactification
$g$ of $g_0$ as in Theorem~3. Here the desingularization is
given by the embedded resolution of $g_0^{-1}(0)\subset(\C^3,0)$
constructed below.

Blow-up the origin of $\C^3$ with $E_0$ the exceptional
divisor. This contains $\mcup_iL_i\cup C$ as its intersection with
the proper transform of $g_0^{-1}(0)$. At each singular point $P_i$ of
$L_i\cup C\subset E_0=\PP^2$, the pull-back of $g_0$ can be written
locally as
$$\bl(v^3(v-u^3)-w^4\br)w^{12},$$
using appropriate analytic local coordinates $u,v,w$. Here $E_0$ is
locally defined by $w=0$, and $u,v$ induce local coordinates of $E_0$
such that $C$ and $L_i$ are respectively defined by $v=0$ and
$v=u^3$ locally on $E_0$.
(Note that we have $w^4$ in the above function since $\deg h''=16$.
The following argument about the point-center blow-ups does not work
well unless $\deg h''=16$.)
We repeat point-center blow-ups three times at the singular point $P_i$.
Here $u,v,w$ are respectively substituted by $u$, $uv$, $uw$ each time.
After these three blow-ups, we get
$$\bl(v^3(v-1)-w^4\br)u^{48}w^{12}.$$
Here the proper transform $E'_0$ of $E_0$ is locally defined by $w=0$,
and the proper transforms of $C$, $L_i$, which will be denoted
respectively by $C'$, $L'_i$, are defined by $v=0$ and $v=1$ locally on
$E'_0$. So $C'$ and $L'_i$ do not intersect each other.
Let $E_i$ denote the exceptional divisor of the last blow-up for each
$i=1,2,3$. This is locally defined by $u=0$ using the above coordinates
after taking the three blow-ups, and transversally intersects $C'$ and
$L'_i$ as is seen by the above description.

The total transform of $g_0^{-1}(0)$ has still singularities along $C'$.
These can be resolved by repeating the blow-ups with center isomorphic to
$C'$ four times.
Indeed, the pull-back of of $g_0$ is generically given by the function
$$(v^3-w^4)w^{12},$$
after restricting to a hyperplane transversal to $C'$.
Here $v,w$ are respectively replaced with $vw,w$ by the first blow-ups,
and by $v,vw$ by the remaining three blow-ups. We do not have a problem
at the intersection point of $C'$ and $E_i$, since the intersection is
transversal as is seen by the above equation. However, the calculation
at the intersection of $C'$ with the proper transform of $h''{}^{-1}(0)$
is rather non-trivial. (The latter does not intersect $E_i$ for $i=1,2,3$
by the assumption on $h''$.) Using appropriate analytic local coordinates
$u,v,w$, the pull-back of $g_0$ can be written as
$$(v^3-uw^4)w^{12},$$
where $E'_0$, $C'$, and the intersection of $E'_0$ with the proper
transform of $h''{}^{-1}(0)$ are respectively defined by $w=0$, $v=w=0$,
and $u=w=0$. By the successive blow-ups, $u,v,w$ are substituted by
$u,vw,w$ or $u,v,vw$ depending on the two affine charts each time.
By the first blow-up, we get
$$(v^3-uw)w^{15}\q\h{and}\q(1-uvw^4)v^{15}w^{12},$$
on the two affine charts. Here we do not have to consider the second, since
$1-uvw^4\ne 0$ if $w=0$. By the second blow-up, we then get
$$(v^3w^2-u)w^{16}\q\h{and}\q(v^2-uw)v^{16}w^{15}.$$
Here we do not have to consider the first, since $(v^3w^2-u)w^{16}$
defines a divisor with normal crossings. The argument is similar
for the third and fourth blow-ups.

Let $E_4$ and $E_5$ respectively denote the exceptional divisor of
the first and the last blow-up of the successive four blow-ups.
Let $D_i$ be the proper transform of $E_i$ in $X$ for $i=0,\dots,5$.
These are the irreducible components with multiplicity divisible by 3,
and $D_4$ does not intersect the irreducible components with
multiplicity non-divisible by 3.
Moreover, $D_{\{i,j\}}:=D_i\cap D_j$ does not intersect the
irreducible components with multiplicity non-divisible by 3 if and
only if $4\in\{i,j\}$ (i.e. $\{i,j\}=\{0,4\}$, $\{4,5\}$,
$\{i,4\}$ with $i=1,2,3$).
Here $D_{\{0,4\}}$ and $D_{\{4,5\}}$ are isomorphic to the original
elliptic curve $C$. We have the unramified coverings of degree 3
$$\Dt_4\to D_4,\q\Dt_{\{0,4\}}\to D_{\{0,4\}},
\q\Dt_{\{4,5\}}\to D_{\{4,5\}},$$
which are compatible with the base changes by the inclusions
$$D_{\{0,4\}}\into D_4,\q D_{\{4,5\}}\into D_4,$$
and also by the canonical projections
$$D_4\to D_{\{0,4\}},\q D_4\to D_{\{4.5\}}.$$
These coverings are {\it non-trivial} by condition~$(A)$.
(Indeed, we apply Proposition~(3.11) to the case where $Z$, $E_i$ and
$h_i$ in Proposition~(3.11) are respectively $D_{\{0,4\}}$, $D_0$ and
$h(u,v,1)w^4$ using the coordinates $u,v,w$ as in Example~(4.2).
Then the non-triviality follows from the remark after condition~$(A)$.)
On the other hand, $D_{\{i,4\}}$ is $\PP^1$ for $i=1,2,3$, and we have
the triviality of the unramified covering
$$\Dt_{\{i,4\}}\to D_{\{i,4\}}\q(i=1,2,3).$$

Setting $b_{f,\la}^j:=\dim B_{f,\la}^j$,
$c_{f,\la}^{\,j}:=\dim C_{f,\la}^j$, we then get for
$\om=\exp(\pm 2\pi i/3)$
$$\aligned b_{f,\om}^0=0,\q b_{f,\om}^1=3,\q b_{f,\om}^2=6,\\
c_{f,\om}^0=1,\q c_{f,\om}^1=5,\q c_{f,\om}^2=6.\endaligned$$
Hence $\chi(B_{f,\om}^{\ssb})\ne\chi(C_{f,\om}^{\ssb})$, and we have
 $\nu_{g_0,\om}^2=\chi(B_{f,\om}^{\ssb})=3$ by Theorem~3.

A similar argument shows that $b_{f,\om}^j=c_{f,\om}^j$ and hence
$\nu_{g_0,\om}^2=2$ in case condition~$(A)$ is {\it not} satisfied,
i.e. if $\al=1$.
This shows that there is no simple formula for $\nu_{g_0,\la}^j$
using only the combinatorial data of the desingularization of $g_0$
in general. So Theorem~1 follows.

\ms\nin
{\bf 4.4.~Example.}
Assume that $f$ is obtained by taking the minimal resolution of a
good projective compactification $f':X'\to\De$ of a germ of a
holomorphic function at $0\in\C^2$ defined by
$$g_0:=(x^{2a}+y^2)(x^2+y^{2a})\q\h{for}\,\,\,a\ge 2.$$
In this case, $f$ is obtained by repeating point-center blow-ups
$2a-1$ times, where all the exceptional divisors have even
multiplicities, but the proper transforms of the irreducible
components of $g^{-1}(0)$ have multiplicity 1.
(This coincides with the resolution obtained by taking a smooth
subdivision of the dual fan of the Newton polygon.)

We have $B_{f,-1}^{\ssb}=C_{f,-1}^{\ssb}$ for $\la=-1$ by Theorem~4,
and moreover
$$\dim C_{f,-1}^0=2a-3,\q \dim C_{f,-1}^1=2a-2.$$
So we get $\nu_{g_0,-1}^1=1$ by Theorem~3.
This assertion also follows from a
theorem in [St2] for the mixed Hodge numbers of the Milnor cohomology
in the non-degenerate Newton boundary case with $\dim X=2$.
(This example shows that the estimate in [MT], which is given by
$\dim C^j_{f,\la}$, is not very good in general.)
Note that some related argument using a $\Q$-resolution is given
in [MM].

\end{document}